\magnification=\magstep1
\input amstex
\documentstyle{amsppt}
\catcode`\@=11 \loadmathfont{rsfs}
\def\mycal{\mathfont@\rsfs}
\csname rsfs \endcsname \catcode`\@=\active

\vsize=6.5in

\topmatter 
\title Asymptotic orthogonalization \\ of subalgebras  in II$_1$ factors  
\endtitle
\author SORIN POPA \endauthor

\rightheadtext{Orthogonalization of subalgebras}

\affil     {\it  University of California, Los Angeles} \endaffil

\address Math.Dept., UCLA, Los Angeles, CA 90095-1555, popa\@math.ucla.edu \endaddress

\thanks Supported by NSF Grant DMS-1700344, a Simons Fellowship and Chaire FSMP-FMJH 2016  \endthanks

\abstract  Let $M$ be a  II$_1$ factor with a von Neumann subalgebra 
$Q\subset M$ that has infinite index under any projection in $Q'\cap M$ (e.g., if $Q'\cap M$ is diffuse; 
or if  $Q$ is an irreducible subfactor with infinite Jones index). 
We prove that 
given any separable subalgebra $B$ of the ultrapower II$_1$ factor $M^\omega$, for a non-principal ultrafilter $\omega$ on $\Bbb N$, 
there exists a unitary element $u\in M^\omega$ such that $uBu^*$ is orthogonal to $Q^\omega$. 
\endabstract

\endtopmatter

\document

\heading 1. Introduction \endheading

We continue in this paper the study  of approximate independence properties for subalgebras in II$_1$ factors from 
[P7,8]. This time, we investigate the possibility of ``orthogonalizing''  two subalgebras of a II$_1$ factor  
via asymptotic unitary conjugacy of one of them, but uniformly with respect to the other.  

Recall in this respect that two $^*$-subalgebras $N_1, N_2$ in a II$_1$ factor $N$ 
are called {\it orthogonal} (as in [P1]), or $1$-{\it independent} (as in [P7]), if  $\tau(x_1x_2)=\tau(x_1)\tau(x_2)$, $\forall x_1\in N_1$, 
$x_2\in N_2$,  $\tau$ denoting the (unique) trace state on the ambient II$_1$ factor.

Thus, given the von Neumann subalgebras $B, Q$ 
of a II$_1$ factor $M$, the problem we are interested in is to find unitary elements $(u_n)_n \subset M$ such that, when viewing $u=(u_n)_n$ as an element  
in the ultrapower II$_1$ factor $M^\omega$ for some non-principal ultrafilter $\omega$ on $\Bbb N$ ([W]), the algebras $uBu^*$ and $Q^\omega$ are 
orthogonal. While this cannot of course be done if  $Q$ is equal to $M$ and $B\neq \Bbb C$, or $Q$ merely ``virtually equal'' to $M$ 
with $\text{\rm dim} (B)$ large enough, 
we will prove that once $Q$ has ``uniform infinite index'' in $M$ and $B$ is separable,  
then such asymptotic orthogonalisation can be obtained. For instance, if $Q$ is an irreducible subfactor of infinite Jones index and $M$ itself 
is separable, then there exists a unitary element $u\in M^\omega$ 
such that $uMu^*  \perp Q^\omega$. 

The {\it uniform infinite index} condition for a subalgebra $Q$ in a II$_1$ factor $M$ 
that we'll require is that for any non-zero projection $p\in Q'\cap M$, 
the [PP]-index of the inclusion $Qp\subset pMp$ be infinite. This condition, which  we have also used in [P9], 
is equivalent to the condition that the centralizer of $M$ in the Jones basic construction algebra $\langle M, Q \rangle$  for $Q\subset M$ 
contains no non-zero finite projections. This amounts to $M\not\prec_M Q$ in the sense of 
``intertwining by bimodules'' terminology (2.4 in [P6]). 

We will in fact investigate this asymptotic orthogonalization problem in the more general case when 
the unitaries $(u_n)_n$ are subject to constraints, being required to lie in some other given von Neumann subalgebra $P\subset M$.   
Thus, we will prove that if $P\subset M$ is any irreducible subfactor such that $P\not\prec_MQ$, then one can indeed find $u\in \Cal U(P^\omega)$ 
such that $uBu^*\perp Q^\omega$. An example of such a situation is when $P, Q$ are irreducible  subfactors of $M$ 
with $[M:P]<\infty$ and $[M:Q]=\infty$.

More generally, we prove the following: 

\proclaim{1.1. Theorem} Let $M_n$ be a sequence of finite factors, with $\dim M_n \rightarrow \infty$. For each $n$, let   
$Q_n \subset M_n$ be a von Neumann subalgebra and $P_n \subset M_n$ be an irreducible subfactor.  
Let $\omega$ be a non-principal ultrafilter on $\Bbb N$.  Denote by $\text{\bf M}$ the ultraproduct 
$\text{\rm II}_1$ factor $\Pi_\omega M_n$ with  $\text{\bf Q}:=\Pi_\omega Q_n$,  $\text{\bf P}:=\Pi_\omega P_n$ 
viewed as von Neumann subalgebras in $\text{\bf M}$.  
Assume $\text{\bf P}\not\prec_{\text{\bf M}} \text{\bf Q}$. Then, given any separable von Neumann subalgebra $B\subset \text{\bf M}$, 
there exists a unitary element $u\in \text{\bf P}$ such that $uBu^*$ is orthogonal to 
$\text{\bf Q}$.  
\endproclaim

For the above condition $\text{\bf P}\not\prec_{\text{\bf M}} \text{\bf Q}$ to be satisfied, it is sufficient that $P_n \not\prec_{M_n} Q_n$, 
$\forall n$, or that  $\text{\bf Q}'\cap \text{\bf P}=\Pi_\omega (Q_n'\cap M_n)$ be diffuse (see Proposition 2.1 below). The condition is also satisfied 
if $M_n$ are II$_1$ factors, with $P_n=M_n$ and $Q_n\subset M_n$ are irreducible subfactors satisfying $\lim_n [M_n:Q_n]=\infty$. 
It is of course satisfied as well when $P_n=M_n$ and $Q_n$ are abelian, $\forall n$. But in fact, as we will show in Remark 2.4, 
this particular case of Theorem 1.1 can be immediately derived from results in [P4]. 

As mentioned before, when applied to the case all $M_n$ are equal to the same II$_1$ factor $M$ and all $Q_n\subset M$ are equal, 
with $P_n=M$, the above theorem gives: 

\proclaim{1.2. Corollary} Let $M$ be a $\text{\rm II}_1$ factor, $Q\subset M$ a von Neumann subalgebra such that $M\not\prec_M Q$ 
and $B\subset M$ a separable von Neumann subalgebra. Then there exists a unitary element $u\in M^\omega$ such 
that $uBu^*\perp Q^\omega$. 
\endproclaim 

The above result shows in particular that once a countably generated II$_1$ factor can 
be embedded in the ultrapower $R^\omega$ of the hyperfinite II$_1$ factor $R$, then one can actually embed it so that to be 
orthogonal to the ultraproduct of an arbitrary sequence of  irreducible subfactors $Q_n\subset R$ with $\lim_n [R: Q_n] =\infty$. 
This fact may be of interest in relation to  Connes Approximate Embedding conjecture.

Since orthogonality (or $1$-independence) between subalgebras   is the first stage of  $n$-independence, it is 
natural to push the above statement even further, trying to  find 
the unitary $u\in M^\omega$ so that $uBu^*$ becomes  $n\geq 2$ independent (or even free independent) 
to $Q^\omega$. This interesting problem remains open for now. 

Questions about ``rotating'' via unitary conjugacy a  subalgebra $B$ in a II$_1$ factor $M$ so that to become (approximately) 
orthogonal to another subalgebra $Q\subset M$ have been first considered in [P1]. The case  
when the algebra $B$ is $2$-dimensional (so ``smallest possible'') has been studied in [P3], notably in the case 
$Q$ is a subfactor of finite Jones index. We will comment on this and other related problems in the last section of this paper, 
where more applications to Theorem 1.1 will be mentioned as well. 

\vskip.05in

{\it Acknowledgement}. This work  has been completed during my stay at the D\'eparte- ment de Math\'ematique de l'Universit\'e de Paris Sud 
and Institut de Math\'ematique de Jussieu during 
the academic year 2016-2017. I want to thank   C. Houdayer,  G. Pisier,
G. Skandalis and S. Vassout for their kind hospitality and support. I am also grateful to F. Le Maitre for his  
interest in orthogonalization problems,  which renewed my own interest in this topic and motivated me  
to write this paper.   

\heading 2. Proof of the main results  \endheading

For notations and terminology used hereafter, we send the reader to ([P7,8]; also [AP]  
for basics in II$_1$ factors theory).

We begin by proving some of the criteria for the condition 
$\text{\bf P}\not\prec_{\text{\bf M}} \text{\bf Q}$ to hold true, that we mentioned in 
the introduction. We will be under the same general assumptions and notations as in Theorem 1.1, 
which are recalled for convenience.

\proclaim{2.1. Proposition} Let $M_n$ be a sequence of finite factors, with $\dim M_n \rightarrow \infty$. For each $n$, let   
$Q_n \subset M_n$ be a von Neumann subalgebra and $P_n \subset M_n$ be an irreducible subfactor.  
Let $\omega$ be a non-principal ultrafilter on $\Bbb N$ and denote by $\text{\bf M}$ the ultraproduct 
$\text{\rm II}_1$ factor $\Pi_\omega M_n$, with  $\text{\bf Q}:=\Pi_\omega Q_n$,  $\text{\bf P}:=\Pi_\omega P_n$ 
viewed as von Neumann subalgebras in $\text{\bf M}$.  Assume one of the following conditions is satisfied: 
\vskip.05in

$1^\circ$ $P_n \not\prec_{M_n} Q_n$, $\forall n$.  

$2^\circ$ $\Pi_\omega (Q_n'\cap M_n)$ is diffuse $($E.g.: all $Q_n$ abelian; or all $M_n$ are finite dimensional factors 
with $P_n=M_n$ and $Q_n$ subfactors with $\lim_n(\dim M_n/\dim Q_n)=\infty)$. 

$3^\circ$ $M_n$ are $\text{\rm II}_1$ factors, $P_n$ is equal to $M_n$ and $Q_n\subset M_n$ are irreducible subfactors with $\lim_n [M_n:Q_n]=\infty$. 

\vskip.05in
\noindent
Then $\text{\bf P}\not\prec_{\text{\bf M}} \text{\bf Q}$. 

\endproclaim

\noindent
{\bf 2.2. Remark.}  Another criterion  for the condition $\text{\bf P}\not\prec_{\text{\bf M}} \text{\bf Q}$ 
to hold true is that $P_n, Q_n\subset M_n$ be irreducible II$_1$ 
subfactors of finite Jones index satisfying $\lim_n [M_n:Q_n]/[M_n:P_n]=\infty$. This result, whose proof requires a lengthier analysis,   
will be discussed in a forthcoming  paper,  which will in fact address a variety of intertwining problems. 

\vskip.1in

\noindent
{\it Proof of Proposition 2.1}.  $1^\circ$ Let $x_1, ..., x_m \in (\text{\bf M})_1$ and $\varepsilon >0$. Let $x_{i,n}\in (M_n)_1$ be so that 
$x_i=(x_{i,n})_n$, $1\leq i \leq m$. Since $P_n\not\prec_{M_n} Q_n$, by (Theorem 2.1 in [P6]) there 
exists a unitary element $u_n\in P_n$ such that $\|E_{Q_n}(x^*_{i,n}u_nx_{j,n})\|_2\leq 2^{-n}$, $1\leq i,j \leq m$. 
Thus, if we let $u=(u_n)_n$ then $u$ is a unitary element in $\Pi_\omega P_n=\text{\bf P}$ satisfying $E_{\text{\bf Q}}(x_i^*ux_j)=0$, 
$\forall i,j$. By (Theorem 2.1), this shows that $\text{\bf P}\not\prec_{\text{\bf M}} \text{\bf Q}$. 

$2^\circ$ By (Theorem 2.1 in [P6]), the condition $\text{\bf P}\prec_{\text{\bf M}} \text{\bf Q}$ would imply 
that there exists an intertwining partial isometry $v\in \text{\bf M}$ between $\text{\bf P}$ and $\text{\bf Q}$. 
Since $\text{\bf P}\subset \text{\bf M}$ is irreducible and $\text{\bf Q}'\cap \text{\bf M}=\Pi_\omega (Q_n'\cap M_n)$ is diffuse, 
this implies $v^*v\in \text{\bf P}$ and $v\text{\bf P}v^* = Q_0 q'$ for some $\text{\it wo}$-closed subalgebra $Q_0\subset \text{\bf Q}$ and 
some projection 
$q'\in \text{\bf Q}'\cap \text{\bf M}$, with $vv^*=1_{Q_0}q'$. 
But the relative commutant $(Q_0q')'\cap vv^* \text{\bf M}vv^*$ contains 
$vv^*(\text{\bf Q}'\cap \text{\bf M})vv^*$ and   is thus diffuse, while by spatiality $v\text{\bf P}v^*$ 
has trivial relative commutant in $vv^*\text{\bf M}vv^*$, a contradiction. 

$3^\circ$ Since $\lim_n [M_n:Q_n]=\infty$, we have $[\text{\bf M}: \text{\bf Q}]=\infty$ (see e.g. [PP]). 
Since we also have $\text{\bf Q}'\cap \text{\bf M}=\Bbb C$, this implies $\text{\bf M}'\cap \langle \text{\bf M},  \text{\bf Q}\rangle = \Bbb C$. But   
$\langle \text{\bf M},  \text{\bf Q}\rangle$ is type II$_\infty$, so the only non-zero projection in $\text{\bf M}'\cap \langle \text{\bf M},  \text{\bf Q}\rangle$ 
is $1_{\langle \text{\bf M},  \text{\bf Q}\rangle}$, which is not finite in $\langle \text{\bf M},  \text{\bf Q}\rangle$. 

\hfill 
$\square$ 

\vskip.05in 

To prove Theorem 1.1, we first show that for any 
$F\subset \text{\bf M} \ominus \Bbb C$ finite and any $\varepsilon > 0$ there exists a unitary element $v\in \text{\bf P}$ 
such that the expectation onto $\text{\bf Q}$ of any element in $vFv^*$  is 
$\varepsilon$-close to $0$ in the Hilbert norm given by the trace. Such unitary $v$  will be constructed by patching together  ``incremental pieces'' of it, 
along the lines of the technique developed in [P5,7,8]. The theorem then follows by a ``diagonalisation along $\omega$'' procedure 
of this local result, as in ([P7], [P8]).

\proclaim{2.3. Lemma} Let $M$ be a $\text{\rm II}_1$ factor, 
$Q\subset M$ a von Neumann subalgebra and $P \subset M$ an irreducible subfactor such that $P\not\prec_M Q$. 
Given any finite set $F=F^*\subset (M\ominus \Bbb C1)_1$ and any $\varepsilon_0 >0$, there 
exists a unitary element $v_0\in P$ such that $\|E_Q(v_0xv_0^*)\|_2^2 \leq \varepsilon_0$, for all $x\in F$. 
\endproclaim

\noindent
{\it Proof}.  Let $\omega$ be a non-principal ultrafilter on $\Bbb N$ 
and denote by $\Cal M=\langle M^\omega, Q^\omega \rangle$ the semi finite von Neumann algebra 
associated with Jones basic construction for $Q^\omega \subset M^\omega$. 
Thus, $\Cal M = J{Q^\omega}' J = \overline{\text{\rm sp}M^\omega e M^\omega}^w \subset \Cal B(L^2M^\omega)$, 
where $e=e_{Q^\omega} \in \Cal B(L^2M^\omega)$ denotes the orthogonal projection of $L^2M^\omega$ onto $L^2Q^\omega$ 
and $J$ is the canonical canonical conjugation on $L^2M^\omega$. 

Recall that the projection $e$ satisfies the condition $eye=E_{Q^\omega}(y)$ for any $y\in M^\omega$. 
Recall also that the semi-finite 
von Neumann algebra $\Cal M$ is endowed with a canonical normal faithful semi-finite trace $Tr$, 
satisfying the condition $Tr(xey)=\tau(xy)$, for all $x, y\in M^\omega$.   

Fix $\varepsilon >0$ such that $\varepsilon < \varepsilon_0$. Denote by ${\mycal W}$ the set of partial isometries $v\in P^\omega$ 
with the property that $vv^*=v^*v$ and which 
satisfy the conditions: 
$$
 \|E_{Q^\omega}(vxv^*)\|_2^2 \leq \varepsilon \tau(v^*v), \ \ \tau(vv^*x) = 0, \tag 1
$$
for all $x\in F$. We endow ${\mycal W}$ with the order $\leq$  
in which $v_1 \leq v_2$ if $v_1=v_2v_1^*v_1$. 

$({\mycal W}, \leq)$ is then clearly inductively ordered and we let $v\in {\mycal W}$ be a maximal element. 
Assume $\tau(v^*v) < 1$ and denote $p = 1 - v^*v$.  Notice right away that $\tau(pFp)=0$.   

Let $w$ be a
partial isometry in $pP^\omega p$ with $w^*w=ww^*$ and denote $u=v+w$. Then $u$ is a partial isometry in $P^\omega$ with $u^*u=uu^*$. 
We will show that one can make an appropriate choice 
$w\neq 0$ such that $u = v + w$ lies in ${\mycal W}$. 
We will construct  such a $w$ by first choosing its support $q=ww^*=w^*w$, then choosing the ``phase $w$'' above $q$. 

Note first that by writing  
$eux^*u^*euxu^*$ as $e(v+w)x^*(v+w)^*e(v+w)x(v+w)^*$ and developing into the sum of 16 terms, we get 
$$
\|E_{Q^\omega}(uxu^*)\|_2^2=Tr(eux^*u^*euxu^*)  \tag 2
$$
$$
=Tr(evx^*v^*evxv^*) + \Sigma_1+\Sigma_2+\Sigma_3+ \Sigma_4, 
$$
where $\Sigma_i$ denotes the sum of terms having $i$ appearances of elements from $\{w, w^*\}$, $1\leq i \leq 4$. 
Thus, there are four terms in $\Sigma_1$, six in $\Sigma_2$, four in $\Sigma_3$, and one in $\Sigma_4$. 

Let us first take care of the terms $Tr(X)$ with $X$ containing a pattern of the form $... ewxw^*e...$, or $...ewx^*w^*e...$, for a given $x\in F$. There are 
seven such terms: the one in $\Sigma_4$, all four in $\Sigma_3$, and two in $\Sigma_2$.  We denote by $\Sigma'$ the sum 
of the absolute values of these terms. Note that for each such $X$, we have $|Tr(X)|=|Tr(wxw^*ey)|$ for some $y\in (M^\omega)_1$ 
and thus, by applying the Cauchy-Schwartz inequality and taking into account the definition of $Tr$, we get the estimate 
$$
|Tr(X)|=|Tr(....ewxw^*e..)| =|Tr(wxw^*ey)| \tag 3
$$
$$
\leq (Tr(ewx^*w^*wxw^*e))^{1/2}(Tr(qey^*yeq))^{1/2} \leq \|qxq\|_2\|q\|_2, 
$$
where the last inequality is due to the fact that $Tr(qey^*yeq)\leq Tr(qeq)=\tau(q)$ 
and $Tr(ewx^*w^*wxw^*e)=Tr(ewx^*qxw^*e)=\tau(wx^*qxw)=\tau(qx^*qxq)$. 

By  (Corollary $2.2.(i)$ in [P4]), the irreducible subfactor  $pP^\omega p$ of  the II$_1$ factor $pM^\omega p$ 
contains a diffuse abelian subalgebra that's $2$-independent to $pFp$  
with respect to the trace state $\tau( \cdot )/\tau(p)$  on $pM^\omega p$. This implies that  
there exists a projection $q\in \Cal P(pP^\omega p)$ of trace $\tau(q)=\varepsilon^2 \tau(p)^2/64$ such that $\tau(qx)=0$ 
and $\|qxq\|_2^2/\tau(p)=(\tau(q)/\tau(p))^2\tau(x^*x)/\tau(p)$, for all $x\in pFp$ and thus for all $x\in F$ as well (because $q\leq p$). 
Thus, for each such $x\in F$ one has $\|qxq\|^2_2 = (\varepsilon^4\tau(p)^2/64^2) \tau(x^*x)\leq\varepsilon^2 \tau(q)/64$.  

It follows that $\|qxq\|_2\leq \varepsilon \tau(q)^{1/2}/8$, $\forall x\in F$. Hence,  
for this choice of $q$, the  right hand side term in $(3)$ will be majorized by $\varepsilon \tau(q)/8$.  
By summing up over the seven terms in $\Sigma'$, we get $\Sigma' \leq 7\varepsilon \tau(q)/8$. 

We now estimate the sum $\Sigma''$ of  $|Tr(X)|$ with $X$ running over the remaining four terms in $\Sigma_2$, and the 
sum $\Sigma'''$ of four terms $|Tr(X)|$ with $X$ having only one occurrence of $w, w^*$ (i.e., the sum of the 
absolute values of the terms in $\Sigma_1$). We'll show that one can choose 
the ``phase $w$'' above the (fixed by now) projection $q$ in a way that makes $\Sigma''+\Sigma'''$ be majorized by $\varepsilon\tau(q)/16$. 

At this point, it is convenient to enumerate the elements in $F=\{x_1, ..., x_n\}$. For each $i=1, 2, ..., n$ we have  
$$
\Sigma''=|Tr(ewx_i^*v^*evx_iw^*)|+|Tr(evx_i^*w^*ewx_iv^*)| \tag 4 
$$
$$
+|Tr(ewx_i^*v^*ewx_iv^*)|+|Tr(evx_i^*w^*evx_iw^*)|
$$
$$
=|Tr(w^*ewY_{1,i})|+|Tr(w^*ewY_{2,i})|
$$
$$
+|Tr(wY_{3,i}wY_{4,i})|+|Tr(w^*Y_{5,i}w^*Y_{6,i})|
$$
where each one of the terms $Y_{j,i}$ depends on $x_i\in F$ and 
belongs to the set  $S_0$ $:=q((M^\omega)_1e(M^\omega)_1)q$ $\subset qL^2(\Cal M, Tr)q$. 

Note that, as $i=1, 2, ..., n$,  the number 
of possible indices $(j,i)$ in  $(4)$ is $6n$. Note also that there are $2n$ terms of the form $|Tr(w^*ewY)|$, 
$n$ terms of the form $|Tr(wXwY)|$ and $n$ terms of the form $|Tr(w^*Xw^*Y)|$, which by 
using the fact that $|Tr(w^*Xw^*Y)|=|Tr(wX^*wY^*)|$ we can view as $n$ additional terms of the form 
$|Tr(wXwY)|$. In all this, the elements $X, Y$ belong to $S_0\subset qL^2(\Cal M, Tr)q$, and are thus 
bounded in operator norm by $1$ and are supported (from left and right) by projections of trace $Tr$ majorized by $1$.  

Similarly, as $i$ runs over $\{1, 2, ..., n\}$, the four terms in $\Sigma'''$ give rise to $4n$ terms of the form $|Tr(wX)|$, for some  
$X\in S_0$. Note that by the definition of $Tr$, each one of these terms is equal to $|\tau(wy)|$ for some $y\in (qM^\omega q)_1$. 

We want to prove that for any $\delta>0$ there exists $w\in \Cal U(qP^\omega q)$ such that each one of the 
above $8n$ terms is less than $\delta$. 

To take care of the terms in $\Sigma'''$, note that by results in ([P4] or [P8]) for any 
finite subset $E\subset qM^\omega q$, 
there exists a finite dimensional subfactor $P_0\subset qP^\omega q$ such that $\|E_{P_0'\cap qM^\omega q}(y)-\tau(y)/\tau(q) q\|_2\leq 
\delta \tau(q)/2$, $\forall y\in E$. By applying this to the elements in $\Sigma'''$, which are of the form $|\tau(wy)|$ 
with $y$ running over a certain finite set $E\subset (qM^\omega q)_1$, and using the Cauchy-Schwartz inequality, one obtains 
that for each unitary element  $w\in N:=P_0'\cap qP^\omega q$ of trace satisfying 
$|\tau(w)|\leq \delta \tau(q)/2$, we have 
$$
|\tau(wy)| = |\tau(E_{P_0'\cap qM^\omega q}(wy))| = |\tau(wE_{P_0'\cap qM^\omega q}(y))| 
$$
$$
\leq  |\tau(w(E_{P_0'\cap qM^\omega q}(y) -(\tau(y)/\tau(q))q)|+ |\tau(w)| |\tau(y)|
$$
$$
\leq \delta \tau(q)/2 + \delta \tau(q)/2 = \delta \tau(q), 
$$
forall $y\in E$. Taking $\delta$ sufficiently small, one obtains that  for any $1\leq i \leq n$ one has 
$\Sigma''' \leq \varepsilon \tau(q)/32$, for any unitary element $w\in N$ 
satisfying $|\tau(w)|\leq \delta\tau(q)/2$. 

Finally, let us take care of the terms in $\Sigma''$. To do this, recall that we are under  the assumption 
$P\not\prec_M Q$,  which in turn implies $P^\omega\not\prec_{M^\omega} Q^\omega$.  
Thus, $(P^\omega)' \cap \Cal M$ contains no finite non-zero projections of $\Cal M=\langle M^\omega, Q^\omega \rangle$ and so  
$N'\cap q\Cal M q$ contains no finite non-zero projections of  $\Cal M$ either. 

To estimate the terms in $\Sigma''$, we first show that for any $\delta_0>0$ and any two $m$-tuples of elements $(Z_1, ..., Z_m)$, 
$(Z_1', ..., Z_m')$ in $S_0 \cap \Cal M_+$, 
there exists a unitary element $w\in N$ such that 
$$
\Sigma_i Tr(w^*Z_iwZ_i')\leq \delta_0. \tag 5
$$

To see this, let $\Cal H$ denote the Hilbert space $L^2(q\Cal Mq, Tr)^{\oplus m}$ and note that we have a unitary 
representation $\Cal U(N)\ni w \mapsto \pi(w) \in \Cal U(\Cal H)$, which on an $m$-tuple 
$X=(X_i)_{i=1}^{m}\in \Cal H$ acts by $\pi(w)(X)=(w^*X_iw)_i$.  

Now note that this representation has no  
(non-zero) fixed point. Indeed, for if $X\in \Cal H$ satisfies $\pi(w)(X)=X$, $\forall w\in \Cal U(N)$, 
then on each component $X_i\in L^2(q\Cal Mq, Tr)$ of $X$ we would have $w^*X_iw=X_i$, $\forall w$. 
Thus $X_iw=wX_i$ and since the unitaries of $N$ span linearly the algebra $N$, 
this would imply $X_i\in N'\cap L^2(q\Cal Mq, Tr)$. Hence, $X_i^*X_i\in N'\cap L^1(q\Cal Mq, Tr)$ and therefore all spectral 
projections of $X_i^*X_i$ corresponding to intervals $[t, \infty)$ with $t>0$ would be projections 
of finite trace in $N'\cap q\Cal Mq$, forcing them all to be equal to $0$. Thus, $X_i=0$ for all $i$. 

With this in mind, denote by $K_Z\subset \Cal H$ the weak closure of the convex hull of the set $\{\pi(w)(Z) \mid w\in \Cal U(N)\}$, 
where $Z=(Z_1, ..., Z_m)$ is viewed as an element in $\Cal H$. Since $K_Z$ is bounded and weakly closed, it is weakly compact, 
so it has a unique element $Z^0\in K_Z$ of minimal norm $\| \ \|_{2,Tr}$. Since $K_Z$ is invariant to $\pi(w)$ and $\|\pi(w)(Z^0)\|_{2,Tr}=\|Z^0\|_{2,Tr}$, 
it follows that $\pi(w)(Z^0)=Z^0$. But we have shown that $\pi$ has no non-zero fixed points, and so $0=Z^0\in K_Z$. 

Let us deduce from this that if $Z=(Z_i)_i, Z'=(Z'_i)_i$ are the two $m$-tuples of positive elements in $S_0$, then 
we can find $w\in \Cal U(N)$ such that $(5)$ holds true. Indeed, for if there would exist $\delta_0>0$ such that $\Sigma_i Tr(\pi(w)(Z_i)Z_i')\geq \delta_0$, 
$\forall w\in \Cal U(N)$, then by taking convex combinations and weak closure, one would get $0=\langle Z^0, Z'\rangle \geq \delta_0$,  
a contradiction. 

This finishes the proof of $(5)$. Note that by taking for one of the $i$  the elements $Y_i, Y_i'$ to be equal to $e$, 
one can get $w\in \Cal U(N)$ to also satisfy  $|\tau(w)|^2\leq \delta_0$.  

We will now use this fact to prove that, given any $m$-tuples $(X_i)_i, (Y_i)_i, (X'_i)_i,$ $(Y'_i)_i \in S_0^m$ (not necessarily 
having positive operators as entries), there exists $w\in \Cal U(N)$ such that $|Tr(w^*X_iwX'_i)|\leq \delta_0$, 
$|Tr(wY_iwY_i')|\leq \delta_0$, for all $i$. Indeed, because if we denote by $e_i$ the left support of $X'_i$ and 
$f_i$ the left support of $Y'_i$, then by the Cauchy-Schwartz inequality we simultaneously have for all $i$ the estimates 
$$
|Tr(w^*X_iwX'_i)|^2\leq Tr(w^*X_i^*X_iwX'_i{X'_i}^*)Tr(e_i) 
\leq Tr(w^*X_i^*X_iwX'_i{X'_i}^*), 
$$
and respectively 
$$
|Tr(wY_iwY_i')|^2\leq Tr(w^*Y_i^*Y_iwY'_i{Y'_i}^*)Tr(f_i)
 \leq Tr(w^*Y_i^*Y_iwY'_i{Y'_i}^*). 
$$

Since all $X_i^*X_i, X'_i{X'_i}^*, Y_i^*Y_i, Y'_i{Y'_i}^*$ are positive elements in $S_0$, we can now apply $(5)$ 
to deduce that  there exist $w\in \Cal U(N)$ of arbitrarily small trace such that all the $4n$ elements 
appearing in $\Sigma''$ for $i=1, 2, ..., n$ are arbitrarily small, making $\Sigma'' \leq \varepsilon \tau(q)/32$, $\forall i$. 
 
Altogether, we then get  for $u=v+w$ the estimate 
$$
\|E_{Q^\omega}(uxu^*)\|_2^2 = \|E_{Q^\omega}(vxv^*)\|_2^2+\Sigma_1 + \Sigma_2+\Sigma_3+\Sigma_4 
$$
$$
\leq \varepsilon \tau(vv^*) + \Sigma'+\Sigma''+\Sigma''' 
$$
$$
\leq \varepsilon \tau(vv^*)+15\varepsilon\tau(ww^*)/16+ 2  \varepsilon\tau(ww^*)/32=\varepsilon\tau(uu^*), 
$$
for any $x\in F$. Thus, $u\in {\mycal W}$ and $u\geq v$, $u\neq v$, contradicting the maximality of $v$. 

This shows that $v$ must be a unitary element. Thus, if we represent $v\in P^\omega$ as a sequence of unitary 
elements $(v_n)_n$ in $P$, then we have 
$$
\lim_{n\rightarrow \omega} \|E_Q(v_nxv_n^*)\|_2^2 = \|E_Q(vxv^*)\|_2^2 \leq \varepsilon < \varepsilon_0, 
$$
for all $x\in F$. Thus, if we let $v_0=v_n$  for some large enough $n$, then $v_0$ satisfies 
$\|E_Q(v_0xv_0^*)\|_2^2 \leq \varepsilon_0$, for all $x\in F$.

\hfill 
$\square$ 

\vskip.1in 

\noindent
{\it Proof of Theorem 1.1}. Let $\{b_j\}_j\subset B$ be a sequence of elements 
that's dense in $(B)_1$ in the Hilbert norm $\| \  \|_2$.  By applying the above lemma to the factor $\text{\bf M}=\Pi_\omega M_n$, 
with its von Neumann subalgebra $\text{\bf Q}=\Pi_\omega Q_n$ and its 
irreducible subfactor $\text{\bf P}=\Pi_\omega P_n$, the finite set $F=\{b_1, ..., b_m\}$ and 
$\varepsilon = 2^{-m-1}$, one gets a unitary element $u_m\in \text{\bf P}$ such that $\|E_{\text{\bf Q}}(u_mb_ju_m^*)\|_2^2\leq 2^{-m-1}$ 
for all $1\leq j \leq m$. 

Let us now take representations $b_j=(b_{j,n})_n \in B \subset \text{\bf M}$ and $u_m=(u_{m,n})_n$, with $u_{m,n}\in \Cal U(P_n)$ 
and $b_{j,n}\in (M_n)_1$. Thus, we have 

$$
\lim_{n \rightarrow \omega}\|E_{Q_n}(u_{m,n}b_{j,n}u_{m,n}^*)\|_2^2  \leq 2^{-m-1}.
$$ 

Denote by $V_m$ the set of all $n\in \Bbb N$ such that $\|E_{Q_n}(u_{m,n}b_{j,n}u_{m,n}^*)\|_2^2  < 2^{-m}$, 
for all $1\leq j \leq m$. Note that $V_m$ corresponds to an open closed neighborhood of $\omega$ in $\Omega$, 
under the identification $\ell^\infty\Bbb N = C(\Omega)$. Let now $W_m$, $m\geq 0$, be defined recursively as follows: 
$W_0=\Bbb N$ and $W_{m+1}=W_n \cap V_{v+1}\cap \{n\in \Bbb N \mid n > \min W_n\}$. 
Note that, with the same identification as before, $W_m$ is a strictly decreasing sequence of neighborhoods of $\omega$ in $\Omega$. 

Define $u=(u_m)_m$ by letting $u_n=u_{m, n}$ for $n\in W_{m-1}\setminus W_m$. It is then easy to see 
that $u$ is a unitary element in $\text{\bf P}$ and that $\lim_\omega \|E_{Q_n}(u_nb_{j,n}u_n^*)\|_2^2=0$ 
for any $j$. In other words, $E_{\text{\bf Q}}(ub_ju^*)=0$, $\forall j$. By the density of the set $\{b_j\}_j$ in $(B)_1$, it follows 
that $uBu^*\perp \text{\bf Q}$. 

\hfill 
$\square$ 

\vskip.05in
\noindent
{\it Proof of Corollary 1.2}. This is just the case where all $M_n$ are equal to the same II$_1$ factor $M$, all $P_n\subset M$ are equal to $M$ 
and $Q_n=Q$, $\forall n$, of Theorem 1.1.  

\hfill 
$\square$

\vskip.05in
\noindent
{\bf 2.4. Remark}. Let us note here that the case $P_n=M_n$ and $Q_n\subset M_n$ abelian, $\forall n$, 
in Theorem 1.1 can be easily deduced directly from 
the main Theorem in [P4]. To see this, note first that it is sufficient to 
prove the statement for a larger $\text{\bf Q}$, and since we can embed each $Q_n$ in a 
MASA (maximal abelian $^*$-subalgebra) of $M_n$, it follows that it is sufficient to 
settle the case when $\text{\bf Q}=\Pi_\omega Q_n$ has diffuse center $Z$.  
Then note that by [P4] there exists a Haar unitary $v \in \text{\bf M}$ 
that's free independent to $B$. Since in an ultrapower II$_1$ factor any two Haar unitaries 
are unitary conjugate, since $Z$ was assumed diffuse, 
it follows that there exists $u\in \Cal U(\text{\bf M})$ such that $A_0:=u\{v\}''u^*\subset  Z =Z(\text{\rm Q})$. 
Thus, $uBu^*$ is free independent to $A_0\subset \text{\bf Q}$. But then, by Kesten's Theorem [K], 
if $\{p_k\}_k$ is any finite partition of $1$ with projections in $A_0$ of trace $\leq \varepsilon^2/2$, then 
one has as in [P7] the estimate $\|\Sigma_k p_k x p_k\|\leq \varepsilon$, $\forall x\in (uBu^*)_1$ 
with $\tau(x)=0$. Thus, since $[A_0, \text{\bf Q}]=0$, we have $E_{\text{\bf Q}}(x)=\Sigma_k p_k E_{\text{\bf Q}}(x)p_k=
E_{\text{\bf Q}}(\Sigma_k p_kxp_k)$ has norm $\leq \varepsilon$, with $\varepsilon>0$ arbitrary, implying that 
$x\perp \text{\bf Q}$, i.e., $uBu^*\perp \text{\bf Q}$.

\heading 3. Further comments  \endheading

\noindent
{\it 3.1. Initial work on orthogonal subalgebras}. The orthogonalization relation between subalgebras $B, Q$ of a II$_1$ factor $M$, as well as 
questions about conjugating a subalgebra $B$ by a unitary element $u\in M$ such that $uBu^*$ becomes orthogonal to $Q$, 
have been first considered in [P1]. They were used in that paper as a tool for calculating normalizers, and more generally the intertwining space 
between subalgebras, in the spirit of what has later become the {\it intertwining by bimodules techniques} [P6]. 
For instance, it was shown in [P1] that if a unitary $u\in M$ has the property that $uA_0u^*\perp Q$ for some 
diffuse abelian subalgebra $A_0\subset Q$, then $u$ is perpendicular to the normalizer of $Q$ in $M$. This allowed 
to  prove several indecomposability properties (e.g., absence of  Cartan subalgebras) for ultraproduct II$_1$ factors and for  
free group factors with uncountable  number of  generators. 

Related to asymptotic orthogonalization,   
it has been shown in (Lemma 2.5 in [P1] and Corollary 2.4 in [P9]) that  
in order for a unitary $u$ to be orthogonal 
to the intertwining space $\Cal I(P, Q)$ between subalgebras $P, Q \subset M$, it is necessary and sufficient 
that there exists a diffuse subalgebra $B\subset Q^\omega$ such that $uBu^* \perp P^\omega$.

\vskip.05in

\noindent
{\it 3.2. The orthogonalization problem}.  Given von Neumann subalgebras $B, Q$ in a II$_1$ factor $M$, 
the problem of finding a unitary $u\in M$ such that $uBu^*\perp Q$ 
will be called the {\it orthogonalization  problem} for $B, Q\subset M$. 
For a given von Neumann subalgebra $Q\subset M$ and a fixed finite dimensional subalgebra $B$, with its trace inherited from $M$,  
all isomorphic copies of $B$ are unitary conjugate in $M$. Thus,  the orthogonalization problem becomes a question about whether there exist 
copies of $B$ that are perpendicular to $Q$. This provides a source of isomorphism invariants for the inclusion $Q\subset M$. 

\vskip.05in
\noindent
{\it 3.3. The two dimensional case}. 
The simplest case of this problem is when $B$ is two dimensional, 
i.e., $B=\Bbb Cf+\Bbb C(1-f)$, 
for some projection $f\in M$ of trace $\tau(f)=\alpha$. Thus, since all projections of same trace  
are unitary conjugate in $M$, the question of whether $B$ can be unitarily conjugated to an algebra 
orthogonal to the von Neumann subalgebra $Q\subset M$ becomes: for what $\alpha \in (0,1]$ does there exist $f\in \Cal P(M)$ such 
that $E_Q(f)=\alpha 1$.  

This problem   
has been systematically investigated in [P3], where the answers   
are formulated in terms of the invariant 
$$
\Lambda(Q\subset M):=\{\alpha \in [0,1] \mid \exists f\in \Cal P(M), E_Q(f)=\alpha 1\}.$$ 
Similarly, one denotes its approximate version  by 
$$
\Lambda_{app}(Q\subset M):=\{\alpha \in ([0,1] \mid \forall \varepsilon >0, \exists f\in \Cal P(M), \|E_Q(f)-\alpha 1\|_2 \leq \varepsilon\}. 
$$   
It obviously coincides with $\Lambda(Q^\omega \subset M^\omega)$ and satisfies the property: $\alpha\in \Lambda_{app}$  
if and only if the algebra $B=\Bbb Cf+\Bbb C(1-f)$, with $f$ a projection of trace $\alpha$,  
can be asymptotically conjugated to an algebra orthogonal to $Q$. 
Note also that $\Lambda_{app}(Q^\omega \subset M^\omega)=\Lambda(Q^\omega\subset M^\omega)$. 

It was already noticed in [P3] that if $Q\subset M$ is an irreducible subfactor of infinite index, 
then (2.4 in [PP]) can be used to show that $\Lambda_{app}(Q\subset M)=\Lambda(Q^\omega\subset M^\omega)=[0,1]$. 
This is of course implied by Theorem 1.1, which in fact applies to all cases when $M\not\prec_M Q$, for instance 
when $Q=A$ is a MASA in $M$. 

The calculation of the invariant $\Lambda(Q\subset M)$ is in general quite difficult, but some partial answers could be obtained in [P3] 
in several particular cases. For instance, if $Q=A$ is a MASA then $\Lambda(A\subset M)$ contains all rationals in $[0,1]$ and if  
in addition $\Cal N_M(A)''$ is of type II$_1$ then $\Lambda(A\subset M)=[0,1]$ (this ought to be the case for any MASA). 

In the finite index case, the results obtained in [P3] depend on weather $[M:Q]<4$ (thus $[M:Q]=4\cos^2 (\pi/n+2)$ 
for some $n\geq 1$, by Jones celebrated results in [J]), 
or $[M:Q]\geq 4$. To describe them,  
denote $[M:Q]^{-1}=\lambda$  and define recursively the polynomials $P_k(x)$   
by $P_{-1}\equiv 1$, $P_0\equiv 1$, $P_{k+1}(x)=P_k(x)-xP_{k-1}(x)$, $k\geq 0$. 
Then, (Theorem in [P3]) shows that if $[M:Q]=\lambda^{-1}=4\cos^2 (\pi/n+2)$   then 
we have $\Lambda(Q\subset M)=\{0\}\cup \{\lambda P_{k-1}/\lambda P_k(\lambda) \mid 0\leq k \leq n-1\}$. 
While if $[M:Q]\geq 4$ and we let $0<t\leq 1/2$ be so that $t(1-t)=\lambda$, then $\Lambda(Q\subset M)\cap (0,t) = 
\{\lambda P_{k-1}(\lambda)/P_k(\lambda) \mid k \geq 0\}$.

Thus, the situation is quite rigid when the index is under the threshold $4$, with the set $\Lambda$ being finite and completely 
understood. While above the threshold $4$ 
the set is always infinite, being completely determined in the intervals $[0, t) \cup (1-t, 1]$, but with the calculation of 
$\Lambda(Q\subset M)\cap [t, 1-t]$ still open in general. It is interesting to note that in both cases (Theorem in [P3]) 
provides the following uniqueness result as well:   any two projections $f_1, f_2\in M$ 
satisfying $E_{Q}(f_1)=\alpha 1=E_Q(f_2)$, with $\alpha=\lambda P_{k-1}(\lambda)/ P_k(\lambda)$ for some $k\geq 0$, are conjugate 
by a unitary element in $Q$. Moreover, since $[M^\omega: Q^\omega]=[M:Q]$, the results in [P3] show that 
$\Lambda(Q\subset M)=\Lambda_{app}(Q\subset M)$ for index $<4$ and $\Lambda(Q\subset M) \cap (0,t)=\Lambda_{app}(Q\subset M)\cap (0,t)$ 
for index $\geq 4$, with any projection that's close to expect on a scalar in these sets being close to a projection that 
actually expects on a scalar. 

The results in [P3] also show that for subfactors $Q\subset M$ of index $\geq 4$ that are {\it locally trivial}, 
i.e., for which $Q'\cap M=\Bbb C f + \Bbb C(1-f)$, $fMf=Qf$, $(1-f)M(1-f)=Q(1-f)$, with $\tau(f)=t\leq 1/2$ where $t(1-t)=\lambda=[M:Q]^{-1}$, 
the invariant $\Lambda(Q\subset M)$ contains no points in the interval $(t, 1-t)$, being equal to 
the set $\{0, t\} \cup \{\lambda P_{k-1}(\lambda)/P_k(\lambda) \mid k \geq 0\}$ when intersected with $[0,1/2]$. This is in particular the case 
when $Q=\{e_n \mid n\geq 1\}''\subset \{e_n \mid n\geq 0\}''=M$ is a subfactor generated by Jones projections of trace $\tau(e_n)=\lambda<1/4$. 

The opposite phenomenon may hold true for subfactors $Q\subset M$ with graph $A_\infty$ and $[M:Q]=\lambda^{-1}>4$. Namely, it is quite 
possible that in all such cases 
one has  $[t, 1-t]\subset \Lambda(Q\subset M)$. As a supporting evidence, 
consider the standard representation $\Cal Q^{st}\subset^{\Cal E} \Cal M^{st}$ of $Q\subset M$, as described in [P10]. This is an inclusion 
of discrete type I von Neumann algebras (i.e., direct sums of type I$_\infty$ factors) 
with inclusion graph $A_\infty$ and a conditional expectation $\Cal E$ having the property that 
$\Cal E_{|M}=E^M_Q$. Also, $\Cal E$ is the unique expectation preserving the trace $Tr$ on $\Cal M^{st}$ whose weights 
are proportional to the square roots of indices of irreducible subfactors 
in the Jones tower. Since $Q\subset M$ is embedded with commuting squares into $\Cal Q^{st}\subset^{\Cal E} \Cal M^{st}$, 
one has $\Lambda(Q\subset M)\subset \Lambda(\Cal Q^{st} \subset^{\Cal E} \Cal M^{st})$, and it is an easy exercise to see that 
the latter contains the entire interval $[t, 1-t]$ (however, it is not clear how one could ``push down'' into $M$ 
a projection $p\in \Cal M^{st}$ that satisfies $\Cal E(p)=s1$ with $s\in [t, 1-t]$).

\vskip.05in
\noindent
{\it 3.4. The finite dimensional case}.  The orthogonalization problem  is certainly interesting for any finite dimensional abelian 
algebra $B=\Sigma_{i=1}^n \Bbb Cf_i$, with $f_1, ..., f_n$ a partition of $1$ with projections in $M$ 
of trace $\tau(f_i)=\alpha_i$, beyond the case $n=2$.  
To state this problem properly, for each $n\geq 2$ we consider the set $F_n$ of all 
$n$-tuples $\alpha=(\alpha_1, ..., \alpha_n)$ with $0\leq \alpha_i \leq 1$ and $\Sigma_i \alpha_i=1$.  
If $Q$ is a von Neumann subalgebra of $M$, we denote 
$$
\Lambda_{n-1}(Q\subset M):=\{(\alpha_i)_i\in F_n \mid 
\exists f_1, ..., f_n \in \Cal P(M), \Sigma_i f_i =1, E_Q(f_i)=\alpha_i, \forall i\}.
$$  
Thus, $\Lambda_1(Q\subset M)=\{(\beta, 1-\beta) \mid \beta \in \Lambda(Q\subset M)\}$. 
Also, any $(\alpha_i)_i\in \Lambda_{n-1}(Q\subset M)$ produces elements in $\Lambda_{k-1}(Q\subset M)$ with $k\leq n-2$ 
by taking $k$-tuples $(\beta_j)_{j=1}^k$  corresponding to partitions of $\{1, ..., n\}$ into $k$ subsets $S_j$ 
and letting $\beta_j=\Sigma_{i\in S_j} \alpha_i$. Thus, the restrictions on $\Lambda=\Lambda_1$ propagate into 
a set of restrictions for $\Lambda_n$, $n\geq 2$. In particular, any entry of an element in $\Lambda_n(Q\subset M)$ 
for $n\geq 1$ is at least $[M:Q]^{-1}$. 

The question here is to calculate (or at least estimate) the invariants $\Lambda_n(Q\subset M)$ for all $n\geq 1$. 
The case when $Q$ is an irreducible 
subfactor of finite index is particularly interesting. 
One source of  $(n+1)$-tuples $\alpha\in \Lambda_n(Q\subset M)$ is to take irreducible subfactors $P\subset Q$ and look 
for partitions of $1$ with $n+1$ projections $f_1, ..., f_n, f_{n+1}$ in $P'\cap M$. Indeed, because then $E^M_Q(f_i)\in P'\cap Q=\Bbb C1$ 
must be scalars. Such $P\subset Q$ can be taken to be an irreducible subfactor obtained by reducing inclusions from 
a Jones tunnel $Q_{-m} \subset Q_{-m+1} \subset ... \subset Q_{-1}\subset Q \subset M$ associated with $Q\subset M$ by a minimal projection 
in $Q_{-m}'\cap Q$. Thus, an $(n+1)$-tuple in $\Lambda_n$ will appear whenever one has an $(n+1)$-point 
in the principal  graphs of $Q\subset M$, $M\subset \langle M, Q \rangle$. For instance, any triple point in $\Gamma_{Q\subset M}$, 
will produce an element in $\Lambda_2(Q\subset M)$, whose entries are proportional to square roots of indices 
of the corresponding irreducible subfactors in the Jones tunnel/tower. When combined with the restrictions on entries coming from 
the obstructions on $\Lambda(Q\subset M)=\Lambda_1(Q\subset M)$, this can provide restrictions on the existence of triple points (and thus 
of graphs of subfactors). More generally, one can apply this same reasoning to the universal graph of $\Gamma^u_{Q\subset M}$, as defined in 
[P10]. 

Another interesting question 
for subfactors of finite index $Q\subset M$ is whether $\Lambda(Q\subset M)$ (more generally $\Lambda_n(Q\subset M)$, $n\geq 1$) depends 
solely on the standard invariant $\Cal G_{Q\subset M}$ of the subfactor $Q\subset M$. The results in [P3] 
show that if $[M:Q]<4$, then in fact $\Lambda(Q\subset M)=\Lambda_1(Q\subset M)$  only depends 
on $[M:Q]$. This may be the case for all $\Lambda_n(Q\subset M)$, when $[M:Q]\leq 4$, but 
it is quite unclear for index $>4$, where however the irreducibly condition $Q'\cap M=\Bbb C$ should be imposed. 
A test case is when $\Gamma_{Q\subset M}=A_\infty$  (i.e., when $Q\subset M$ has TLJ standard invariant) with the index running over the 
interval $(4, \infty)$ (cf. [P5]).

\vskip.05in
\noindent
{\it 3.5.  Unitaries in orthonormal basis}.  One case of particular interest is to decide whether for some given $2\leq k \leq n$ one can have 
$\alpha=(\alpha_i)_i\in \Lambda_n(Q\subset M)$ with $\alpha_1= \alpha_2=...=\alpha_k=\lambda=[M:Q]^{-1}$. Thus, in such a case one has 
$s\lambda \in \Lambda(Q\subset M)$ for any $s=0, 1, 2, ..., k$. 
By (Proposition 1.7 in [PP]), this is equivalent to whether $Q_{-1} \subset Q$ has an orthonormal basis $\{m_i\}_i$  
with the first $k$ terms $1=m_1, ..., m_k$ being unitary elements. 

If $[M:Q] \not\in \Bbb N$ and $n$ denotes its integer part then,  
as pointed out in (1.4.2$^\circ$ of [PP]), the formula $\Sigma_i m_im_i^*=[M:Q] 1$ implies $k\leq n-1$. 

The [P3] restrictions on $\Lambda(Q\subset M)$ can be used to obtain further restrictions 
on the maximal number of unitaries that can appear in an orthonormal basis of $M$ over $Q$. 
Indeed, if the index $[M:Q]$ is less than $4$ but $\neq 3$, then the equations $\lambda P_{m-1}(\lambda)/P_m(\lambda)=2\lambda$ 
do not have solutions for $\lambda^{-1}=4\cos^2 (\pi/n+2)$ and $n \neq 2, 4$. If in turn $[M:Q]>4$ 
and we let $[M:Q]=n+\varepsilon$, with $1>\varepsilon >0$, then having $n-1$ mutually orthogonal projections $f_1, ..., f_{n-1}\in M$ 
with $E_Q(f_i)=\lambda=[M:Q]^{-1}$ would imply that $f=1-\Sigma_i f_i$ satisfies $\tau(f)=1-(n-1)\lambda  \in \Lambda(Q\subset M)$ and $\tau(f)\geq 
\lambda/(1-\lambda)$. This in turn implies $\varepsilon \geq \lambda/(1-\lambda)$. Thus,  $\varepsilon < \lambda/(1-\lambda)$ 
forces $k\leq n-2$. In particular, this shows that if $4< [M:Q] < 4+ \frac{\sqrt{13}-3}{2}\approx 4.3$, then $k\leq 2$, i.e., there exists at most one unitary 
$u\in M$ with $E_Q(u)=0$. 

In turn, it is an open problem whether an irreducible subfactor $Q\subset M$ with integer index $n\geq 5$ always has an 
orthonormal basis with $n$ unitaries. Note that by a result in [P2], if a subfactor $Q$ of a II$_1$ factor $M$ contains a Cartan subalgebra of $M$, then 
$Q\subset M$ does have an orthonormal basis of unitaries (even if its index is infinite). However, if a subfactor $Q\subset M$ has $A_\infty$-graph 
(e.g., when $Q\subset M$ is constructed by the universal amalgamated free product method in [P5]), 
then finding even one single unitary $u\in M$ with $E_Q(u)=0$ is an open question (see also the conjecture at the end of 3.3 above). 

\vskip.05in
\noindent
{\it 3.6.  A related dilation problem}. The question about whether a scalar $\alpha \in (0,1)$ 
is the expected value on $Q$ of a projection in $M$ is viewed in [P3] as a dilation problem. 
More generally, one can ask this same question for an arbitrary element $b\in Q$ with $0\leq b \leq 1$:  
can $b$ be dilated to a projection $p\in M$, i.e.,  
does there exist a projection  $p\in M$ such 
that $E_Q(p)=b$ ?

Alternatively, one can attempt the calculation of the entire set $E_Q(\Cal P(M))$, or of the set 
$\tilde{\Lambda}(Q\subset M)$ of all functions $g:[0,1]\rightarrow [0,1]$ that can be spectral distributions 
of elements in $E_Q(\Cal P(M))$, i.e., with the property
that there exists $p\in \Cal P(M)$ with $g(\alpha)=\tau(e_\alpha(E_Q(p))$, $\forall \alpha\in [0,1]$. 
These sets should be calculable for subfactors of index $<4$. 

On the other hand, note that if $Q\subset M$ satisfies the infinite index condition $M\not\prec_M Q$, then Corollary 1.2 
easily yields a calculation of the approximate versions of these sets, showing that 
$E_{Q^\omega}(\Cal P(M^\omega))= \{b\in Q^\omega \mid 0\leq b \leq 1\}$ and 
thus $\tilde{\Lambda}(Q^\omega \subset M^\omega)$ is equal to the set of non-increasing functions from $[0,1]$ to $[0,1]$. 
To see this, note first that $E_{Q^\omega}(\Cal P(M^\omega))$ is $\| \ \|_2$-closed (by the usual $\omega$-diagonalisation procedure).  
Then note that any $b\in Q^\omega$ with $0\leq b \leq 1$ can be approximated uniformly by elements of finite spectrum 
$b'=\Sigma_i \alpha_i q_i$ with $q_i\in \Cal P(Q^\omega)$, $0\leq \alpha_i \leq 1$. Finally, as noticed in 3.3 above, any $\alpha_i q_i$ can 
be dilated to a projection $p_i\in q_iM^\omega q_i$, thus $p=\Sigma_i p_i$ satisfies $E_{Q^\omega}(p)=b'$. 

\vskip.3in

\head  References \endhead

\item{[AP]} C. Anantharaman, S. Popa: ``An introduction to II$_1$ factors'', \newline www.math.ucla.edu/$\sim$popa/Books/

\item{[C]} A. Connes: {\it Classification of injective factors},
Ann. of Math., {\bf 104} (1976), 73-115.

\item{[J]} V.F.R. Jones: {\it Index for subfactors}, Invent. Math. {\bf 72} (1983), 1-25. 

\item{[K]} H. Kesten: {\it Symmetric random walks on groups}, Trans. Amer. Math. Soc. {\bf 92} (1959), 336-354.

\item{[PP]} M. Pimsner, S. Popa, {\it Entropy and index for
subfactors}, Annales Scient. Ecole Norm. Sup., {\bf 19} (1986),
57-106.

\item{[P1]} S. Popa, {\it Orthogonal pairs of *-subalgebras in
finite von Neumann algebras}, J. Operator Theory, {\bf 9} (1983),
253-268.

\item{[P2]} S. Popa, {\it Notes on Cartan subalgebras in type
$II_1$ factors}. Mathematica Scandinavica, {\bf 57} (1985),
171-188.

\item{[P3]} S. Popa: {\it Relative dimensions, towers of
projections and commuting squares of subfactors}, Pacific J.
Math., {\bf 137} (1989), 181-207,

\item{[P4]} S. Popa: {\it Free independent sequences in type} II$_1$ {\it factors
and related problems}, Asterisque, {\bf 232} (1995), 187-202.

\item{[P5]} S. Popa, {\it An axiomatization of the lattice of
higher relative commutants of a subfactor}, Invent. Math., {\bf
120} (1995), 427-445.

\item{[P6]} S. Popa: {\it Strong Rigidity of }  II$_1$ {\it Factors
Arising from Malleable Actions of $w$-Rigid Groups} I, Invent. Math.,
{\bf 165} (2006), 369-408.

\item{[P7]} S. Popa: {\it A} II$_1$ {\it factor approach to the Kadison-Singer problem}, 
Comm. Math. Physics. {\bf 332} (2014), 379-414 (math.OA/1303.1424).

\item{[P8]}  S. Popa: {\it Independence properties in subalgebras of ultraproduct} II$_1$ {\it factors}, Journal of Functional Analysis 
{\bf 266} (2014), 5818-5846 (math.OA/1308.3982)

\item{[P9]}  S. Popa:  {\it Constructing MASAs with prescribed properties}, to appear in Kyoto J. of Math, math.OA/1610.08945

\item{[P10]} S. Popa: {\it Some properties of the symmetric enveloping algebras
with applications to amenability and property T},
Documenta Mathematica, {\bf 4} (1999), 665-744.

\item{[V]} D. Voiculescu:  {\it Symmetries of some reduced free product 
$C^*$-algebras},  In: ``Operator algebras and their connections with topology
and ergodic theory'', Lect. Notes in Math. Vol. {\bf 1132}, 566-588 (1985).

\item{[W]} F. B. Wright:  {\it A reduction for algebras of finite type},  Ann. of Math. {\bf 60} (1954), 560 - 570.

\enddocument